\documentclass{amsart}
\usepackage{amsmath}
\usepackage{amssymb}
\usepackage{amsthm}

\DeclareMathOperator{\Vol}{Vol}

\newcommand{\oP}{\overline{P}}

\newcommand{\hg}{\widehat{g}}

\newcommand{\hP}{\widehat{P}}
\newcommand{\hQ}{\widehat{Q}}

\newcommand{\htheta}{\widehat{\theta}}

\newcommand{\lv}{\lvert}
\newcommand{\rv}{\rvert}

\newcommand{\mP}{\mathcal{P}}

\newcommand{\bC}{\mathbb{C}}

\newcommand{\bR}{\mathbb{R}}

\def\sideremark#1{\ifvmode\leavevmode\fi\vadjust{\vbox to0pt{\vss
 \hbox to 0pt{\hskip\hsize\hskip1em
 \vbox{\hsize3cm\tiny\raggedright\pretolerance10000
 \noindent #1\hfill}\hss}\vbox to8pt{\vfil}\vss}}}

\newcommand{\comment}[1]{}

\newtheorem{thm}{Theorem}[section]

\newtheorem{cor}[thm]{Corollary}

\theoremstyle{definition}

\theoremstyle{remark}

\numberwithin{equation}{section}

\begin{document}

\title{Extremal metrics for the ${Q}^\prime$-curvature in three dimensions}
\author{Jeffrey S.\ Case}
  \address{Department of Mathematics, McAllister Building, The Pennsylvania State University, University Park, PA 16802}
  \email{jscase@psu.edu}
\author{Chin-Yu Hsiao}
  \address{Institute of Mathematics, Academia Sinica, 6F, Astronomy-Mathematics Building, No.1, Sec.4, Roosevelt Road, Taipei 10617, Taiwan}
  \thanks{CYH was supported by Taiwan Ministry of Science of Technology project 
103-2115-M-001-001, 104-2628-M-001-003-MY2 and the Golden-Jade fellowship of Kenda Foundation}
 \email{chsiao@math.sinica.edu.tw}
\author{Paul Yang}
  \thanks{PY was partially supported by NSF Grant DMS-1509505}
  \address{Department of Mathematics \\ Princeton University \\ Princeton, NJ 08544}
  \email{yang@math.princeton.edu}
\begin{abstract}
 We construct contact forms with constant $Q^\prime$-curvature on compact three-dimensional CR manifolds which admit a pseudo-Einstein contact form and satisfy some natural positivity conditions.  These contact forms are obtained by minimizing the CR analogue of the $II$-functional from conformal geometry.  Two crucial steps are to show that the $P^\prime$-operator can be regarded as an elliptic pseudodifferential operator and to compute the leading order terms of the asymptotic expansion of the Green's function for $\sqrt{P^\prime}$.
\end{abstract}
\maketitle

\section{Introduction}
\label{sec:intro}

On an even-dimensional manifold $(M^{2n},g)$, the pair $(P,Q)$ of the (critical) GJMS operator $P$ and the (critical) $Q$-curvature $Q$ possesses many of the same properties of the pair $(-\Delta,K)$ on surfaces, where $K$ is the Gauss curvature.  For example, $P$ is a conformally covariant formally self-adjoint operator with leading order term $(-\Delta)^{n/2}$ which annihilates constants~\cite{GJMS1992,GrahamZworski2003} and $Q$ is a Riemannian invariant with leading order term $c_n(-\Delta)^{\frac{n-2}{2}}R$, where $R$ is the scalar curvature, which transforms in a particularly simple way within a conformal class~\cite{Branson1995}: if $\hg=e^{2u}g$, then
\[ e^{n\sigma}\hQ = Q + Pu . \]
In particular, $\int Q$ is conformally invariant on closed even-dimensional manifolds; indeed, it computes the Euler characteristic modulo integrals of pointwise conformal invariants~\cite{Alexakis2012}.  It also follows that metrics of constant $Q$-curvature within a conformal class are in one-to-one correspondence with critical points of the functional
\[ II[u] = \int_M u\,Pu + 2\int_M Qu - \frac{2}{n}\left(\int_M Q\right)\log\left(\frac{1}{\Vol(M)}\int_M e^{nu}\right) . \]
This functional can always be minimized on the two-sphere~\cite{OsgoodPhillipsSarnak1988} and on four-manifolds with positive Yamabe constant and nonnegative Paneitz operator~\cite{Beckner1993,ChangYang1995,Gursky1999}, with important applications to logarithmic functional determinants~\cite{BransonChangYang1992,OsgoodPhillipsSarnak1988} and sharp Onofri-type inequalities~\cite{Beckner1993}.  Due to the parallels between conformal and CR geometry, it is interesting to determine whether a similar pair exists in the latter setting.

Work of Graham and Lee~\cite{GrahamLee1988} and Hirachi~\cite{Hirachi1990} identified CR analogues of the Paneitz operator and $Q$-curvature in dimension three.  However, the kernel of the Paneitz operator contains the (generally infinite-dimensional) space $\mP$ of CR pluriharmonic functions and the total $Q$-curvature is always zero.  In particular, an Onofri-type inequality involving the CR Paneitz operator cannot be satisfied.  Branson, Fontana and Morpurgo overcame this latter issue on the CR spheres by introducing a formally self-adjoint operator $P^\prime$ which is CR covariant on CR pluriharmonic functions and in terms of which one has the sharp Onofri-type inequality
\[ \int u\,P^\prime u + 2\int Q^\prime u - \frac{2}{n+1}\left(\int Q^\prime\right)\log\left(\frac{1}{\Vol(S^{2n+1})}\int e^{(n+1)u}\right) \geq 0 \]
for all $u\in W^{n+1,2}\cap\mP$, where $Q^\prime$ is an explicit dimensional constant~\cite{BransonFontanaMorpurgo2007}.  The construction of $P^\prime$ is analogous to the construction of the $Q$-curvature from the GJMS operators by analytic continuation in the dimension.

It was observed by the first- and third-named authors in dimension three~\cite{CaseYang2012} and by Hirachi in general dimension~\cite{Hirachi2013} that one can define the $P^\prime$-operator on general pseudohermitian manifolds $(M^{2n+1},T^{1,0},\theta)$.  Roughly speaking, if $P_{2n+2}^N$ is the CR GJMS operator of order $2n+2$ on a $(2N+1)$-dimensional manifold, one defines $P^\prime$ as the limit of $\frac{2}{(N-n)}P_{2n+2}^N\rv_\mP$ as $N\to n$.  This is made rigorous by explicit computation in dimension three~\cite{CaseYang2012} and via the ambient metric in general dimension~\cite{Hirachi2013}.  Regarded as a map from $\mP$ to $C^\infty(M)/\mP^\perp$, the $P^\prime$-operator is CR covariant: if $\htheta=e^\sigma\theta$, then $e^{(n+1)\sigma}\hP^\prime=P^\prime$.

If $\theta$ is a pseudo-Einstein contact form (cf.\ \cite{CaseYang2012,Hirachi2013,Lee1988}), then the $P^\prime$-operator is formally self-adjoint and annihilates constants.  Note that if $M^{2n+1}$ is the boundary of a domain in $\bC^{n+1}$, then the defining functions constructed by Fefferman~\cite{Fefferman1976} induce pseudo-Einstein contact forms on $M$.  One can construct a pseudohermitian invariant $Q^\prime$ on pseudo-Einstein manifolds by formally considering the limit $\bigl(\frac{2}{N-n}\bigr)^2P_{2n+2}^N(1)$ as $N\to n$; this can be made rigorous by direct computation in dimension three~\cite{CaseYang2012} and via the ambient metric in general dimension~\cite{Hirachi2013}.  Regarded as $C^\infty(M)/\mP^\perp$-valued, the $Q^\prime$-curvature transforms linearly with a change of contact form: if $\htheta=e^\sigma\theta$ is also pseudo-Einstein, then
\begin{equation}
 \label{eqn:qprime_transformation}
 e^{2(n+1)}\hQ^\prime = Q^\prime + P^\prime(\sigma) .
\end{equation}
Since $\htheta$ is pseudo-Einstein if and only if $\sigma\in\mP$~\cite{Hirachi1990,Lee1988}, this makes sense.  It follows from the properties of $P^\prime$ that $\int Q^\prime$ is independent of the choice of pseudo-Einstein contact form.  Direct computation on $S^{2n+1}$ shows that it is a nontrivial invariant; indeed, in dimension three it is a nonzero multiple of the Burns--Epstein invariant~\cite{CaseYang2012}.  In particular, the pair $(P^\prime,Q^\prime)$ on pseudo-Einstein manifolds has the same properties as the pair $(P,Q)$ on Riemannian manifolds.

If $(M^{2n+1},T^{1,0},\theta)$ is a compact pseudo-Einstein manifold, the self-adjointness of $P^\prime$ and~\eqref{eqn:qprime_transformation} imply that critical points of the functional $II\colon\mP\to\bR$ defined by
\begin{equation}
 \label{eqn:ii_functional}
 II[u] = \int_M u\,P^\prime u + 2\int_M Q^\prime u - \frac{2}{n+1}\left(\int_M Q^\prime\right)\log\left(\frac{1}{\Vol(M)}\int_M e^{(n+1)u}\right)
\end{equation}
are in one-to-one correspondence with pseudo-Einstein contact forms with constant $Q^\prime$-curvature (still regarded as $C^\infty(M)/\mP^\perp$-valued).  The existence and classification of minimizers of the $II$-functional on the standard CR spheres was given by Branson, Fontana and Morpurgo~\cite{BransonFontanaMorpurgo2007}.  In this note, we discuss the main ideas used by the authors to give criteria which guarantee that minimizers exist for the $II$-functional on a given pseudo-Einstein three-manifold~\cite{CaseHsiaoYang2014}.

\begin{thm}
 \label{thm:main_thm}
 Let $(M^3,T^{1,0},\theta)$ be a compact, embeddable pseudo-Einstein three-manifold such that $P^\prime\geq0$ and $\ker P^\prime=\bR$.  Suppose additionally that
 \begin{equation}
  \label{eqn:qprime_bound}
  \int_M Q^\prime\,\theta\wedge d\theta < 16\pi^2 .
 \end{equation}
 Then there exists a function $w\in\mP$ which minimizes the $II$-functional.  Moreover, the contact form $\htheta:=e^w\theta$ is such that $\hQ^\prime$ is constant.
\end{thm}

The assumptions on $P^\prime$ mean that the pairing $(u,v):=\int u\,P^\prime v$ defines a positive definite quadratic form on $\mP$.  It is important to emphasize that the conclusion is that $\hQ_4^\prime$ is constant as a $C^\infty(M)/\mP^\perp$-valued invariant: a local formula for the $Q^\prime$-curvature was given by the first- and third-named authors~\cite{CaseYang2012}, while we observe that, on $S^1\times S^2$ with any of its locally spherical contact structures, there is no pseudo-Einstein contact form with $Q^\prime$ pointwise zero; see~\cite[Section~5]{CaseHsiaoYang2014}.

As in the study of Riemannian four-manifolds (cf.\ \cite{ChangYang1995,Gursky1999}), the hypotheses of Theorem~\ref{thm:main_thm} can be replaced by the nonnegativity of the pseudohermitian scalar curvature and of the CR Paneitz operator.  Indeed, Chanillo, Chiu and the third-named author proved that these assumptions imply that $(M^3,T^{1,0})$ is embeddable~\cite{ChanilloChiuYang2010}; the first- and third-named authors proved that these assumptions imply both that $P^\prime\geq0$ with $\ker P^\prime=\bR$ and that $\int Q^\prime\leq 16\pi^2$ with equality if and only if $(M^3,T^{1,0})$ is CR equivalent to the standard CR three-sphere~\cite{CaseYang2012}; and Branson, Fontana and Morpurgo showed that minimizers of the $II$-functional exist on the standard CR three-sphere~\cite{BransonFontanaMorpurgo2007}.

\begin{cor}
 \label{cor:main_cor}
 Let $(M^3,T^{1,0},\theta)$ be a compact pseudo-Einstein manifold with nonnegative scalar curvature and nonnegative CR Paneitz operator.  Then there exists a function $w\in\mP$ which minimizes the $II$-functional.  Moreover, the contact form $\htheta:=e^w\theta$ is such that $\hQ^\prime$ is constant.
\end{cor}
\section{Sketch of the proof of Theorem~\ref{thm:main_thm}}
\label{sec:sketch}

The proof of Theorem~\ref{thm:main_thm} proceeds analogously to the proof of the corresponding result on four-dimensional Riemannian manifolds~\cite{ChangYang1995} with one important difference: $P^\prime$ is defined as a $C^\infty(M)/\mP^\perp$-valued operator; in particular, it is a nonlocal operator.  Let $\tau\colon C^\infty(M)\to\mP$ be the orthogonal projection with respect to the standard $L^2$-inner product.  A key observation is that the operator $\oP^\prime:=\tau P^\prime \colon \mP\to \mP$ is a self-adjoint elliptic pseudodifferential operator of order $-2$; see~\cite[Theorem~9.1]{CaseHsiaoYang2014}.  This follows from the observation that, while the sublaplacian $\Delta_b$ is subelliptic, the Toeplitz operator $\tau\Delta_b\tau$ is a classical elliptic pseudodifferential operator of order $-1$.  This is achieved by writing $\Delta_b=2\Box_b+iT$, relating $\tau$ to the Szeg{\H o} projector $S$, and using well-known properties of the latter operator (cf.\ \cite{BoutetSjostrand1976,Hsiao2010}).

Since $\int u\,P^\prime v = \int u\,\oP^\prime v$ for all $u,v\in\mP$, it follows that $\oP^\prime$ is a nonnegative operator with $\ker\oP^\prime=\bR$.  In particular, the positive square root $\bigl(\oP^\prime\bigr)^{1/2}$ of $\oP^\prime$ is well-defined and such that $\ker\bigl(\oP^\prime\bigr)^{1/2}=\bR$.  Using the pseudodifferential calculus and the fact that, as a local operator, $P^\prime$ equals $\Delta_b^2$ plus lower order terms~\cite{CaseYang2012}, we then observe that the Green's function of $\bigl(\oP^\prime\bigr)^{1/2}$ is of the form $c\rho^{-2}+O(\rho^{-1-\varepsilon})$ for $\rho^4(z,t)=\lv z\rv^4+t^2$ the Heisenberg pseudo-distance, $\varepsilon\in(0,1)$, and $c$ the same constant as the computation on the three-sphere~\cite{BransonFontanaMorpurgo2007}; for a more precise statement, see~\cite[Theorem~1.3]{CaseHsiaoYang2014}.

From this point, the remaining argument is fairly standard.  The above fact about the Green's function of $\bigl(\oP^\prime\bigr)^{1/2}$ allows us to apply the Adams-type theorem of Fontana and Morpurgo~\cite{FontanaMorpurgo2011} to conclude that the former operator satisfies an Adams-type inequality with the same constant as on the standard CR three-sphere.  This has two important effects.  First, it implies that that $II$-functional is coercive under the additional assumption $\int Q^\prime<16\pi^2$; see~\cite[Theorem~4.1]{CaseHsiaoYang2014}.  Second, it implies that if $w\in W^{2,2}\cap\mP$ satisfies
\[ \tau\left( P^\prime w + Q^\prime - \lambda e^{2w}\right) = 0, \]
then $w\in C^\infty(M)$; see~\cite[Theorem~4.2]{CaseHsiaoYang2014}.  The former assumption allows us to minimize $II$ within $W^{2,2}\cap\mP$ and the latter assumption yields the regularity of the minimizers.  The final conclusion follows from the transformation formula~\eqref{eqn:qprime_transformation} for the $Q^\prime$-curvature.

\subsection*{Acknowledgements}

The authors thank Po-Lam Yung for his careful reading of an early version of the article~\cite{CaseHsiaoYang2014}.  They also thank the Academia Sinica in Taipei and Princeton University for warm hospitality and generous support while this work was being completed.

\bibliographystyle{abbrv}
\bibliography{../bib}

\newcommand{\noopsort}[1]{}
\begin{thebibliography}{10}

\bibitem{Alexakis2012}
S.~Alexakis.
\newblock {\em The decomposition of global conformal invariants}, volume 182 of
  {\em Annals of Mathematics Studies}.
\newblock Princeton University Press, Princeton, NJ, 2012.

\bibitem{Beckner1993}
W.~Beckner.
\newblock Sharp {S}obolev inequalities on the sphere and the
  {M}oser-{T}rudinger inequality.
\newblock {\em Ann. of Math. (2)}, 138(1):213--242, 1993.

\bibitem{BoutetSjostrand1976}
L.~Boutet~de Monvel and J.~Sj{\"o}strand.
\newblock Sur la singularit\'e des noyaux de {B}ergman et de {S}zeg{\H o}.
\newblock In {\em Journ\'ees: \'{E}quations aux {D}\'eriv\'ees {P}artielles de
  {R}ennes (1975)}, pages 123--164. Ast\'erisque, No. 34--35. Soc. Math.
  France, Paris, 1976.

\bibitem{Branson1995}
T.~P. Branson.
\newblock Sharp inequalities, the functional determinant, and the complementary
  series.
\newblock {\em Trans. Amer. Math. Soc.}, 347(10):3671--3742, 1995.

\bibitem{BransonChangYang1992}
T.~P. Branson, S.-Y.~A. Chang, and P.~C. Yang.
\newblock Estimates and extremals for zeta function determinants on
  four-manifolds.
\newblock {\em Comm. Math. Phys.}, 149(2):241--262, 1992.

\bibitem{BransonFontanaMorpurgo2007}
T.~P. Branson, L.~Fontana, and C.~Morpurgo.
\newblock Moser-{T}rudinger and {B}eckner-{O}nofri's inequalities on the {CR}
  sphere.
\newblock {\em Ann. of Math. (2)}, 177(1):1--52, 2013.

\bibitem{CaseHsiaoYang2014}
J.~S. Case, C.-Y. Hsiao, and P.~C. Yang.
\newblock Extremal metrics for the ${Q}^\prime$-curvature in three dimensions.
\newblock \noopsort{2200}Preprint.

\bibitem{CaseYang2012}
J.~S. Case and P.~C. Yang.
\newblock A {P}aneitz-type operator for {CR} pluriharmonic functions.
\newblock {\em Bull. Inst. Math. Acad. Sin. (N.S.)}, 8(3):285--322, 2013.

\bibitem{ChangYang1995}
S.-Y.~A. Chang and P.~C. Yang.
\newblock Extremal metrics of zeta function determinants on {$4$}-manifolds.
\newblock {\em Ann. of Math. (2)}, 142(1):171--212, 1995.

\bibitem{ChanilloChiuYang2010}
S.~Chanillo, H.-L. Chiu, and P.~Yang.
\newblock Embeddability for 3-dimensional {C}auchy-{R}iemann manifolds and {CR}
  {Y}amabe invariants.
\newblock {\em Duke Math. J.}, 161(15):2909--2921, 2012.

\bibitem{Fefferman1976}
C.~Fefferman.
\newblock Monge-{A}mp\`ere equations, the {B}ergman kernel, and geometry of
  pseudoconvex domains.
\newblock {\em Ann. of Math. (2)}, 103(2):395--416, 1976.

\bibitem{FontanaMorpurgo2011}
L.~Fontana and C.~Morpurgo.
\newblock Adams inequalities on measure spaces.
\newblock {\em Adv. Math.}, 226(6):5066--5119, 2011.

\bibitem{GJMS1992}
C.~R. Graham, R.~Jenne, L.~J. Mason, and G.~A.~J. Sparling.
\newblock Conformally invariant powers of the {L}aplacian. {I}. {E}xistence.
\newblock {\em J. London Math. Soc. (2)}, 46(3):557--565, 1992.

\bibitem{GrahamLee1988}
C.~R. Graham and J.~M. Lee.
\newblock Smooth solutions of degenerate {L}aplacians on strictly pseudoconvex
  domains.
\newblock {\em Duke Math. J.}, 57(3):697--720, 1988.

\bibitem{GrahamZworski2003}
C.~R. Graham and M.~Zworski.
\newblock Scattering matrix in conformal geometry.
\newblock {\em Invent. Math.}, 152(1):89--118, 2003.

\bibitem{Gursky1999}
M.~J. Gursky.
\newblock The principal eigenvalue of a conformally invariant differential
  operator, with an application to semilinear elliptic {PDE}.
\newblock {\em Comm. Math. Phys.}, 207(1):131--143, 1999.

\bibitem{Hirachi1990}
K.~Hirachi.
\newblock Scalar pseudo-{H}ermitian invariants and the {S}zeg{\H o} kernel on
  three-dimensional {CR} manifolds.
\newblock In {\em Complex geometry ({O}saka, 1990)}, volume 143 of {\em Lecture
  Notes in Pure and Appl. Math.}, pages 67--76. Dekker, New York, 1993.

\bibitem{Hirachi2013}
K.~Hirachi.
\newblock {$Q$}-prime curvature on {CR} manifolds.
\newblock {\em Differential Geom. Appl.}, 33(suppl.):213--245, 2014.

\bibitem{Hsiao2010}
C.-Y. Hsiao.
\newblock Projections in several complex variables.
\newblock {\em M\'em. Soc. Math. Fr. (N.S.)}, (123):131, 2010.

\bibitem{Lee1988}
J.~M. Lee.
\newblock Pseudo-{E}instein structures on {CR} manifolds.
\newblock {\em Amer. J. Math.}, 110(1):157--178, 1988.

\bibitem{OsgoodPhillipsSarnak1988}
B.~Osgood, R.~Phillips, and P.~Sarnak.
\newblock Extremals of determinants of {L}aplacians.
\newblock {\em J. Funct. Anal.}, 80(1):148--211, 1988.

\end{thebibliography}
\end{document}